\documentclass[12pt]{amsart}

\usepackage{amsmath,amsthm,amssymb,latexsym}
\usepackage{graphicx}
\usepackage{setspace}

\usepackage[all]{xy}
\usepackage{enumerate}
\newtheorem{tmenv}{Theorem}[section]
\newtheorem{crenv}{Corollary}[section]
\newtheorem{lmenv}[crenv]{Lemma}
\newtheorem{exenv}{Example}[section]
\theoremstyle{definition}
\newtheorem{dfenv}{Definition}[section]
\newtheorem{propenv}{Proposition}[section]
\newtheorem*{remenv}{Remark}

\newcommand{\cor}{\begin{crenv}}
	\newcommand{\roc}{\end{crenv}}
\newcommand{\rem}{\begin{remenv}}
	\newcommand{\mer}{\end{remenv}}
\newcommand{\pr}{\begin{proof}}
	\newcommand{\rp}{\end{proof}}
\newcommand{\ex}{\begin{exenv} \rm }
	\newcommand{\xe}{\end{exenv}}
\newcommand{\tm}{\begin{tmenv}}
	\newcommand{\mt}{\end{tmenv}}
\newcommand{\lm}{\begin{lmenv}}
	\newcommand{\ml}{\end{lmenv}}
\newcommand{\prop}{\begin{propenv}}
	\newcommand{\porp}{\end{propenv}}
\newcommand{\af}{\begin{afenv}}
	\newcommand{\fa}{\end{afenv}}
\newcommand{\df}{\begin{dfenv}}
	\newcommand{\fd}{\end{dfenv}}

\begin{document}

\title[A problem equivalent  to counting DAGs on labeled vertices]
{A problem equivalent  to counting directed acyclic graphs on  labeled vertices}

\author{Zs\'ofia Juh\'asz}

\maketitle

\begin{abstract}
	 An encoding of directed acyclic graphs (DAGs) on labeled vertices 
	 is proposed, which is a generalisation of the Prüfer code for labeled trees, if a certain orienation  on the edges of the tree is  introduced. Hence it is shown that the number of sequences	 $S_1, S_2, \ldots, S_{n-1}$ of subsets of $\{1, 2, \ldots, n\}$ with the property that $|\bigcup_{i=1}^kS_i|\leq k$ for every $1\leq k \leq n-1$, is equal to the number of  DAGs on $n$ labeled vertices.  
\end{abstract}

\section{Introduction}

A directed graph $G=(V, E)$ is an ordered pair, where $V$ is the set of vertices and $E\subseteq V\times V$ is the set of (directed) edges. A directed acyclic graph (DAG) is a directed graph that contains no directed cycles. We say that the vertices of a graph are labeled, if every vertex has a distinct element of set $\{1, \ldots , |V|\}$ assigned to, referred to as the label of the vertex. 

In \cite{Robinson1973} R. W. Robinson proved the following recursive formula for the number $a_n$ of DAGs on $n$ labeled  vertices:

\begin{equation} \label{Eqn_Enumerate}
a_n = \sum_{k=1}^n (-1)^{k-1} {n\choose k}2^{k(n-k)} a_{n-k}.
\end{equation}

E. W. Weisstein conjectured \cite{Weisstein} and B. D. McKay showed \cite{McKay} that the number of $n$ by $n$ $(0,1)$-matrices all of whose eigenvalues are positive real numbers is equal to the number of  DAGs on $n$ labeled vertices. 


We introduce another natural  problem, involving counting the number of certain sequences of sets, and show that it is equivalent to the enumeration of DAGs on $n$ labeled vertices. An encoding of DAGs is proposed by a sequence of the out-neighbor sets of the vertices arranged into a certain order. It is shown that this establishes a bijection between the set of  all DAGs on $n$ labeled vertices and the set of all sequences of subsets of $\{1, 2, \ldots , n\}$ satisfying the property that $|\bigcup_{i=1}^kS_i|\leq k$ for every $1\leq k \leq n-1$. Hence the number of latter sequences is equal to the number of  DAGs on $n$ labeled vertices.


\section{Encoding DAGs by sequences of out-neighbor sets}
For any vertices $u$ and $v$ in a directed graph $G=(V, E)$, $v$ is called an out-neighbor of $u$ if $(u, v)\in E$.  

\df \label{Def_Encoding}
Let $G=(V, E)$ be a DAG on the set of labeled vertices $V=\{1, 2, \ldots , n \}$.  Define the \textit{minimal source sequence}  $u_1, u_2, \ldots , u_{n}$ 
of $G$ recursively as follows: let $u_1$ be the source in $G$ with the smallest label. For every $2\leq i\leq n$: let $u_i$ be the source with the smallest label in the graph obtained by deleting vertices $u_1, \ldots , u_{i-1}$ from $G$.  For every $1\leq i\leq n$ denote by $S_{n-i}$ the set of out-neighbors of $u_i$ in $G$. We shall call  $S_1, S_2,  \ldots, S_{n-1}$  the \textit{encoding} of $G$.


\fd

Note that since $S_0=\emptyset$ always holds, $S_0$ does not need to be included in the encoding.

\lm \label{Lm_SourceSeq}
Let  $G=(V, E)$ be a DAG on the set of labeled vertices $V=\{1, 2, \ldots , n \}$ with minimal source sequence $u_1, u_2, \ldots , u_{n}$ and encoding $ S_1, S_2, \ldots, S_{n-1}$. Let  $S_0=\emptyset$. Then $u_1$ is the smallest label vertex in $V\setminus (\bigcup _{i=0}^{n-1}S_i)$ and  for every $2\leq j\leq n$: $u_j$ is the smallest label vertex in $V\setminus (\{u_1, \ldots, u_{j-1}\}\cup (\bigcup _{i=0}^{n-j}S_i))$.
\ml

\pr
By definion of the minimal source sequence and the encoding, $V\setminus (\bigcup _{i=0}^{n-1}S_i)$ is the set of sources in $G$, and for every $2\leq j \leq n$, $V\setminus (\{u_1, \ldots, u_{j-1}\}\cup (\bigcup _{i=0}^{n-j}S_i))$ is exactly the set of sources in the graph obtained by deleting vertices $u_1, \ldots, u_{j-1}$ from $G$. Hence the statement follows.
\rp

\tm \label{Tm_Existence}

A sequence $S_1, \ldots, S_{n-1}$ of subsets of $V=\{1, 2,\ldots , n\}$ is an encoding of a DAG on $n$ labeled vertices if and only if it satisfies the following property:

\begin{equation} \label{Eqn_Property}
|\bigcup_{i=1}^{k}S_i|\leq k \textrm{ for every  } 1\leq k\leq n-1.
\end{equation}
\mt

\pr
$\Rightarrow$
Note that by the definition of the minimal source sequence, for every $1\leq i\leq n-1$: $S_i\subseteq \{u_{n-i+1}, \ldots , u_{n}\}$. Hence for any $1\leq k\leq n-1$:   
$|\bigcup_{i=1}^kS_i|\leq k$.



$\Leftarrow$
Let $S_1, \ldots, S_{n-1}$ be a sequence of subsets of $V$ satisfying Property \ref{Eqn_Property} and let $S_0=\emptyset$. Denote by $w_1$ the smallest label vertex in $V\setminus\bigcup _{i=0}^{n-1}S_i$. For every $2\leq j\leq n$ let $w_{j}$ be the smallest label vertex in $V\setminus (\{w_1, \ldots, w_{j-1}\}\cup (\bigcup _{i=0}^{n-j}S_i))$. By Property \ref{Eqn_Property} $V\setminus (\{w_1, \ldots, w_{j-1}\}\cup(\bigcup _{i=0}^{n-j}S_i))\neq \emptyset$, hence the sequence $w_1, \ldots ,  w_n$ is well-defined, containing each element of $V$ exactly once. Let $G=(V, E)$ the directed graph where for every $1\leq i\leq n$ the set of out-neighbors of $w_i$ is exactly $S_{n-i}$. We prove that $G$ is a DAG and $ S_1, \ldots, S_{n-1}$ is the encoding of $G$. Note that for any edge $(w_p, w_q)$ in $E$, $p<q$. Indeed, if $(w_p, w_q)\in E$ then by definition of $G$, $w_q\in S_{n-p}$ and since by definition,  $w_1$ is the smallest label vertex in $V\setminus (\bigcup _{i=0}^{n-1}S_i)$ and for $1<q\leq n$  $w_q$ is the smallest label vertex in $V\setminus (\{w_1, \ldots, w_{q-1}\}\cup (\bigcup _{i=0}^{n-q}S_i))$, $w_q\notin\bigcup _{i=0}^{n-q}S_i$, hence  $n-p>n-q$ and so $q>p$.  Therefore  $G$ is a DAG. Let $u_1, \ldots ,  u_n$ be the code ordering of the vertices in $G$. By proof by induction we show that for every $1\leq i\leq n$: $u_i=w_i$. By the construction of $G$, the smallest label vertex in  $V\setminus(\bigcup _{i=0}^{n-1}S_i)$ is exactly the source in $G$ with the smallest label, hence $v_1=w_1$. Suppose now that for some $1\leq j\leq n-1$: $u_i=w_i$ holds for every $1\leq i\leq j$. Then $w_{j+1}$ is the smallest label vertex in $V\setminus (\{u_1, \ldots, u_{j}\}\cup(\bigcup _{i=0}^{n-j-1}S_i))$, where $\bigcup _{i=0}^{n-j-1}S_i$ is exactly the set of those vertices in $G$ which are out-neighbors of at least one of $u_{j+1}, \dots , u_n$. Whence $w_{j+1}$ is the smallest label source in the graph obtained by deleting vertices $u_1, \ldots, u_{j}$ from $G$, and so $u_{j+1}=w_{j+1}$. Therefore $u_i=w_i$ and hence $S_i$ is the set of out-neighbors of $u_{n-i}$ for every $1\leq i\leq n$.
\rp

\tm \label{Tm_Uniqueness}
Every sequence $ S_1, \ldots, S_{n-1}$ of subsets of $\{1, 2, \ldots , n\}$ satisfying Property \ref{Eqn_Property} is the encoding of a unique DAG on the set of vertices $\{1, 2, \ldots , n\}$. 
\mt

\pr
By Theorem \ref{Tm_Existence} if $ S_1, \ldots, S_{n-1}$ of subsets of $\{1, 2, \ldots , n\}$ satisfying Property \ref{Eqn_Property} then it is the encoding of some DAG $G=(V, E)$ with $V=\{1, 2, \ldots , n\}$. Suppose it is also the encoding of DAG $G'=(V, E')$. Then by Lemma \ref{Lm_SourceSeq} the minimal source sequences of $G$ and $G'$ are identical. Hence by the definition of the encoding the sets of out-neighbors of each vertex in $V$ are identical in $G$ and $G'$.. Therefore $G=G'$.  
\rp

\cor \label{Cor_Equivalence}
The number of DAGs on $n$ labeled  vertices equals the number of sequences $S_1, \ldots S_{n-1}$  of subsets of $\{1, 2, \ldots , n\}$ satisfying Property \ref{Eqn_Property}.
\roc

\newpage
By Corollary \ref{Cor_Equivalence} and Robinson's recursive formula \ref{Eqn_Enumerate}:

\cor
For the number $a_n$ of sequences $S_1, \ldots , S_{n-1}$  of subsets of $\{1, 2, \ldots , n\}$ with the property that $|\bigcup _{i=1}^k S_i|\leq k$ for every $1\leq k\leq n-1$:

$$ a_n = \sum_{k=1}^n (-1)^{k-1} {n\choose k}2^{k(n-k)} a_{n-k}.$$

\roc


\rem
The encoding of a DAG  on labeled vertices introduced in Definition \ref{Def_Encoding} is a generalisation of the Prüfer code of trees on labeled vertices \cite{Prufer} in the following sense: Given a tree $T=(V, E)$ on $n$ labeled vertices  denote by $a_1,\ldots , a_{n-2}$ the Prüfer sequence of $T$, and for every $1\leq i\leq n-2$ let $v_i$ be the vertex that is removed from $T$ in step $i$ during the creation of the Prüfer code. Denote by $v_{n-1}$ and by $v_n$ the last two remaining vertices in the tree after step $n-2$, where $v_{n-1}$ is the vertex with the smaller label. Denote by $T'$ the DAG obtained by introducing the following orientation on the edges of $T$: for every $1\leq i<j\leq n$ if $v_i$ and $v_j$ are connected by an edge $e$ in $T$ then let $v_i$ be the starting and $v_j$ be the endpoint of $e$, respectively. Then the encoding of $T'$ according to Definition \ref{Def_Encoding} is the sequence $\{a_1\}, \dots , \{a_{n-2}\}, \{v_{n}\}$.
\mer

\bibliographystyle{amsplain}

\end{document}